\documentclass{amsart} 
\usepackage{amscd}
\usepackage{amsmath,empheq}
\usepackage{amsfonts}
\usepackage{amssymb}
\usepackage{mathrsfs}
\usepackage[all]{xy} 
\newtheorem{theorem}{Theorem}

\newtheorem{definition}[theorem]{Definition}
\newcommand{\bb}{\begin{equation}}
\newcommand{\bbb}{\end{equation}}

\newcommand{\calE}{{\mathcal E}}
\newcommand{\calF}{{\mathcal F}}

\newcommand{\calX}{{\mathcal B}}
\newcommand{\calC}{{\mathcal C}}
\newcommand{\intom}{\int_\Omega}
\newcommand{\di}{\displaystyle}

\newcommand{\ri}{\rightarrow}	 
\newenvironment{proof-sketch}{\noindent{\bf Sketch of Proof}\hspace*{1em}}{\qed\bigskip}

\DeclareMathOperator{\divv}{div}
\everymath{\displaystyle}
\newcommand{\RR}{\mathbb R}

\renewcommand{\leq}{\leqslant}

\renewcommand{\geq}{\geqslant}
 
\begin{document} 
\title[Eigenvalue problems with unbalanced growth]{Eigenvalue problems with unbalanced growth: Nonlinear patterns and standing wave solutions} 
\author[I. Farag\'o]{Istv\'an Farag\'o}
\address[I. Farag\'o]{MTA-ELTE NumNet Research Group, Budapest, Hungary \&
Department of Applied Analysis and Computational Mathematics, ELTE University, Budapest, Hungary \&
Department of Differential Equations, Budapest University of Technology and Economics, Budapest, Hungary}
\email{\tt faragois@cs.elte.hu}
\author[D. Repov\v{s}]{Du\v{s}an Repov\v{s}}
\address[D. Repov\v{s}]{Faculty of Education and Faculty of Mathematics and Physics, University of Ljubljana \& Institute of Mathematics, Physics and Mechanics, 1000 Ljubljana, Slovenia}
\email{\tt dusan.repovs@guest.arnes.si}


\begin{abstract}
We consider two classes of nonlinear eigenvalue problems with double-phase energy and lack of compactness. We establish existence and non-existence results and related properties of solutions. Our analysis combines variational methods with the generalized Pucci-Serrin maximum principle.
\end{abstract}
\maketitle

\section{Introduction}
In this paper, motivated by recent advances in the study of nonlinear problems with unbalanced growth,
we are interested in the mathematical analysis of standing wave solutions of some classes of Dirichlet boundary value problems driven by nonhomogeneous differential operators of the type
\bb\label{azor}\divv [\phi'(|\nabla u|^2)\nabla u],\bbb
where $\phi\in C^1(\RR_+,\RR_+)$ has a different growth near zero and at infinity.
Such a behaviour occurs, for instance, if $\phi(t)=2(\sqrt{1+t}-1)$, which corresponds to the prescribed mean curvature differential operator (capillary surface operator), namely
$$\divv\left(\frac{\nabla u}{\sqrt{1+|\nabla u|^2}}\right).$$
More generally, $\phi(t)$ can behave like $t^{q/2}$ for small $t$ and like $t^{p/2}$ for large $t$, where $1<p<q$. Such a growth is fulfilled if
$$\phi(t)=\frac 2p\,[(1+t^{q/2})^{p/q}-1],$$
which generates the differential operator
$$\divv \left((1+|\nabla u|^q)^{(p-q)/q}|\nabla u|^{q-2}\nabla u \right)\,.$$

A case intensively studied in recent years corresponds to
$$\phi(t)=\left\{
\begin{array}{lll}
&\di \frac 2q\,t^{q/2}&\ \mbox{if}\ t<1\\
&\di \frac 2p\,t^{p/2}-\frac{2(q-p)}{pq}&\ \mbox{if}\ t\geq1.
\end{array}\right.$$
It follows that
$$\phi (|\nabla u|^2)\simeq\left\{
\begin{array}{lll}
&\di |\nabla u|^p,&\quad\mbox{if}\ |\nabla u|\gg 1;\\
&\di |\nabla u|^q,&\quad\mbox{if}\ |\nabla u|\ll 1.
\end{array}\right.
$$
This potential produces the $(p,q)$-Laplace operator $\Delta_p+\Delta_q$, which generates a ``double-phase energy" (according to the terminology of Marcellini and Mingione).

We briefly recall in what follows the roots of double-phase problems.
To the best of our knowledge, problems of this type have been first considered by Ball \cite{ball1,ball2} in the context of problems with cavities in nonlinear elasticity.

Let $\Omega\subset\RR^N$ ($N\geq 2$) be a bounded domain  with smooth boundary. If $u:\Omega\to\RR^N$ is the displacement and if $Du$ is the $N\times N$  matrix of the deformation gradient, then the total energy is defined by
\bb\label{paolo}E(u)=\intom f(x,Du(x))dx,\bbb
where  $f=f(x,\xi):\Omega\times\RR^{N\times N}\to\RR$ is quasiconvex with respect to $\xi$. The simplest example considered by Ball is given by functions $f$ of the type
$$f(\xi)=g(\xi)+h({\rm det}\,\xi),$$
where ${\rm det}\,\xi$ is the determinant of the $N\times N$ matrix $\xi$, and $g$, $h$ are nonnegative convex functions, which satisfy the growth conditions
$$g(\xi)\geq c_1\,|\xi|^p;\quad\lim_{t\to+\infty}h(t)=+\infty,$$
where $c_1$ is a positive constant and $1<p<N$. The condition $p< N$ is necessary to study the existence of equilibrium solutions with cavities, that is, minima of the integral \eqref{paolo} that are discontinuous at one point where a cavity forms. In fact, every $u$ with finite energy belongs to the Sobolev space $W^{1,p}(\Omega,\RR^N)$, and thus it is a continuous function if $p>N$. 

In accordance with these problems arising in nonlinear elasticity, Marcellini \cite{marce1,marce2} considered continuous functions $f=f(x,u)$ with  unbalanced growth that satisfy
$$c_1\,|u|^p\leq |f(x,u)|\leq c_2\,(1+|u|^q)\quad\mbox{for all}\ (x,u)\in\Omega\times\RR,$$
where $c_1$, $c_2$ are positive constants and $1\leq p\leq q$.
We also point out the contributions of Baroni, Colombo and Mingione \cite{mingi2,mingi3}
 in the framework of non-autonomous functionals characterized by the fact that the energy density changes its ellipticity and growth properties according to the point.

 These contributions are in relationship with the work of Zhikov \cite{zhikov1}, who described the
behavior of some phenomena arising in nonlinear
elasticity.
In fact, Zhikov intended to provide models for strongly anisotropic materials in the context of homogenisation.
For instance, Zhikov considered the ``double-phase" energy functional defined by
\bb\label{mingfunc}
{\mathcal P}_{p,q}(u) :=\intom (|\nabla u|^p+a(x)|\nabla u|^q)dx,\quad 0\leq a(x)\leq L,\ 1<p<q,
\bbb
where the modulating coefficient $a(x)$ dictates the geometry of the composite made by
two differential materials, with hardening exponents $p$ and $q$, respectively.
The functional ${\mathcal P}_{p,q}$ falls in the realm of the so-called functionals with
nonstandard growth conditions of $(p, q)$--type, according to Marcellini's terminology. These are functionals of the type in \eqref{paolo}, where the energy density satisfies
$$|\xi|^p\leq f(x,\xi)\leq  |\xi|^q+1,\quad 1\leq p\leq q.$$
General models with $(p,q)$-growth in the context of
geometrically constrained problems have been recently studied by De Filippis \cite{cristina}. This seems to be the first work dealing with $(p,q)$-conditions with manifold constraint. Refined regularity results
are proved
in \cite{cristina},
by using an approximation technique relying on estimates obtained through a careful use of difference quotients. Other recent works dealing with nonlinear problems with unbalanced growth (either isotropic or anisotropic) are the papers by  Bahrouni, R\u adulescu and  Repov\v s \cite{bahr}, Cencelj, R\u adulescu and Repov\v s \cite{cencelj}, and
Papageorgiou, R\u adulescu and  Repov\v s \cite{papag}.

The differential operator defined in \eqref{azor} and which is generated by a potential with variable growth was introduced by
Azzollini {\it et al.} \cite{azo1,azo2} in relationship with wide classes of nonlinear PDEs with a variational structure. We refer to Chorfi and R\u adulescu \cite{chorfi} for the study of a related problem driven by this general differential operator. We also refer to the recent monographs \cite{prrbook,radrep} for some of the abstract methods used in the present paper.

\section{Functional setting and main results}
We are concerned with the existence of nontrivial  solutions of the following quasilinear Schr\"odinger problem with double-power nonlinearities:
\bb\label{prelim}
-\divv [\phi'(|\nabla u|^2)\nabla u]+|u|^{\alpha-2}u=\lambda\, |u|^{s-2}u\quad\mbox{in}\ \Omega\subset\RR^N\ (N\geq 2),
\bbb
where $\lambda$ is a positive parameter.

This equation was studied in \cite{azo2} if $\Omega=\RR^N$ and  under the assumption that the reaction dominates the left-hand side of the problem. In fact, Azzollini,  d'Avenia and Pomponio \cite{azo2}
proved that this equation has a nontrivial non-negative radially symmetric solution, provided that $1<p<q<\min\{N,p^*\}$, $1<\alpha\leq p^*q'/p'$, and $\max\{q,\alpha\}<s<p^*$. A crucial tool in their arguments is a certain compactness property of a Sobolev-type space of radially symmetric functions into Lebesgue spaces. Our purpose in this paper is to consider the same equation but on exterior domains of the Euclidean space. We aim to prove
related existence or non-existence results depending on the values of the parameter $\lambda$ and the competition between the left-hand side of Eq.~\eqref{prelim} and its reaction. More precisely, we are first concerned with the following nonlinear eigenvalue problem
\bb\label{problem}\left\{\begin{array}{lll}
&\di -\divv [\phi'(|\nabla u|^2)\nabla u]+|u|^{\alpha-2}u=\lambda\, |u|^{s-2}u&\quad\mbox{in}\ \Omega\\
&\di u=0&\quad\mbox{on}\ \partial\Omega. \end{array}\right.
\bbb
Here, we assume that $\lambda$ is a positive parameter and $\Omega\subset\RR^N$ ($N\geq 2$) is the complement of a bounded domain with smooth boundary.

The existence of solutions of problem \eqref{problem} was studied by Berestycki and Lions \cite{bere} in the case of the Laplace operator and without the presence of the nonlinear term $|u|^{\alpha-2}u$. In this case, the authors assumed that the reaction is a double-power nonlinearity that has a subcritical behaviour at infinity and a supercritical growth near the origin.

In order to describe the main results of this paper, we start with the basic hypotheses.

Throughout this paper we assume that $\alpha$, $p$, $q$ are real numbers satisfying the following hypothesis:
\bb\label{pq}
 1<p<q<N\quad\mbox{and}\quad
 1<\alpha<\frac{p^*q'}{p'}\,.
\bbb

We assume that the  function $\phi:\RR_+\ri\RR_+$ is of class $C^1$ and  has the following properties:

\smallskip
\noindent ($\phi_1$) $\phi (0)=0;$

\smallskip
\noindent ($\phi_2$) there exists $c>0$ such that $\phi(t)\geq ct^{p/2}$ if $t\geq 1$ and   $\phi(t)\geq ct^{q/2}$ if $0\leq t\leq 1$;

\smallskip
\noindent ($\phi_3$) there exists $C>0$ such that $\phi(t)\leq Ct^{p/2}$ if $t\geq 1$ and   $\phi(t)\leq Ct^{q/2}$ if $0\leq t\leq 1$;

\smallskip
\noindent ($\phi_4$) there exists $0<\mu<1$ such that $2t\phi'(t)\leq s\mu\phi (t)$ for all $t\geq 0$;

\smallskip
\noindent ($\phi_5$) the mapping $t\mapsto \phi (t^2)$ is strictly convex.

\smallskip Since our hypotheses allow that $\phi'$ approaches 0, problem \eqref{problem} is  degenerate and no ellipticity condition is assumed.

For all $1\leq r\leq\infty$, we denote by $\|\,\cdot\,\|_r$ the norm  on the Lebesgue space $L^r(\Omega)$.

\begin{definition}\label{d1} Let $\Omega\subseteq\RR^N$ be an open set.
We define the function space $L^p(\Omega)+L^q(\Omega)$ as the completion of $C^\infty_c(\Omega)$ in the norm
$$\|u\|_{L^p+L^q}:=\inf\{\|v\|_p+\|w\|_q;\ v\in L^p(\Omega),\ w\in L^q(\Omega),\ u=v+w\}.$$
\end{definition}
We set $$\|u\|_{p,q}:=\|u\|_{L^p(\Omega)+L^q(\Omega)}.$$

The space $L^p(\Omega)+L^q(\Omega)$ is an Orlicz space and has been intensively studied by Badiale, Pisani and  Rolando \cite[Sect. 2]{badiale}.
This space is a reflexive Banach space, see \cite[Corollary 2.11]{badiale}.
We point out that the space $L^p(\Omega)+L^q(\Omega)$ is of interest only either $p<q$ or $|\Omega|=+\infty$. Indeed, if $p=q$ or $|\Omega|<+\infty$, then $L^q(\Omega)\subseteq L^p(\Omega)$, hence $L^p(\Omega)+L^q(\Omega)=L^p(\Omega)$.

A key role in our arguments is played by the Banach space
$${\mathcal B}:=\overline{C^\infty_c(\Omega)}^{\|\,\cdot\,\|},$$
where $$\|u\|:=\|\nabla u\|_{p,q}+\|u\|_\alpha.$$

As established in Propositions 2.4 and 2.5 of \cite{azo2}, $\calX$ is a reflexive Banach space. Moreover, if $p'<p^*q'$ then for every $1<\alpha\leq p^*q'/p'$, the space $\calX$ is continuously embedded into $L^{p^*}(\Omega)$; see \cite[Theorem 2.6]{azo2} for more details. We point out that the loss of compactness of the Orlicz embeddings in the case of unbounded domains implies refined variational techniques. Some of the papers dealing with problems with lack of compactness on unbounded domains use particular function spaces where the compactness is preserved, such as spaces of radially symmetric functions. Such a situation occurs in \cite{azo2}, where the main existence property is obtained via a compact embedding. We recall that even if the domain is unbounded, standard compact embeddings remain true, for instance if $\Omega$ is ``thin at infinity", in the sense that
$$\lim_{R\ri\infty}\sup\{|\Omega\cap B(x,1)|\,;\ x\in\RR^N,\, |x|=R \}=0.$$
Such a situation does not hold in our case. Indeed, since $\Omega$ is an exterior domain, then it looks like the whole space $\RR^N$ at infinity and, in particular, it is not a thin domain.

\begin{definition}\label{d2}
A solution of problem \eqref{problem} is a function $u\in \calX\setminus\{0\}$ such that for all $v\in\calX$
$$\int_{\Omega}\left(\phi'(|\nabla u|^2)\nabla u\nabla v+|u|^{\alpha-2}uv-\lambda|u|^{s-2}uv\right)dx=0.$$
\end{definition}

The real number $\lambda$ for which problem \eqref{problem} has a nontrivial solution is an eigenvalue and the corresponding $u\in \calX\setminus\{0\}$ is an eigenfunction of the problem. These terms are in accordance with the related notions introduced by Fu\v{c}ik, Ne\v{c}as, Sou\v{c}ek and Sou\v{c}ek \cite[p. 117]{fucik} in the abstract framework of nonlinear operators. Indeed, if we set
$$Su:=
\frac 12\int_{\Omega}\phi (|\nabla u|^2)dx+\frac 1\alpha\int_{\Omega}|u|^\alpha dx\quad\mbox{and}\quad T(u):=\frac{1}{s}\intom |u|^sdx$$
then $\lambda$ is an eigenvalue for the pair $(S,T)$ is and only if there exists a corresponding eigenfunction, namely a solution of problem \eqref{problem} as described by Definition \ref{d2}.

We first prove that problem \eqref{problem} has a solution for any $\lambda>0$, provided that the reaction ``dominates" the growth in the left-hand side. More precisely, we have the following existence result.

\begin{theorem}\label{t1}
Assume that hypotheses \eqref{pq}, ($\phi_1$)-($\phi_5$) are fulfilled, and $\max\{q,\alpha\}<s<p^*$. Then the following properties are true:

(a) problem \eqref{problem} has a nonnegative solution $U$ for all $\lambda>0$;

(b) $U\in C^{1,\mu}(\Omega\cap B_R(0))$ with $\mu =\mu(R)\in(0,1)$;

(c) $U>0$ in $\Omega$.
\end{theorem}

Next, we are concerned with the following nonlinear problem with variable potential and lack of compactness
\bb\label{problem1}\left\{\begin{array}{lll}
&\di -\divv [\phi'(|\nabla u|^2)\nabla u]+|u|^{\alpha-2}u=\lambda\, a(x)\,|u|^{s-2}u&\quad\mbox{in}\ \Omega\\
&\di u=0&\quad\mbox{on}\ \partial\Omega. \end{array}\right.
\bbb

Accordingly, a solution of problem \eqref{problem1} is a function $u\in \calX\setminus\{0\}$ such that for all $v\in\calX$
$$\int_{\Omega}\left(\phi'(|\nabla u|^2)\nabla u\nabla v+|u|^{\alpha-2}uv-\lambda\, a(x)\,|u|^{s-2}uv\right)dx=0.$$

Hypothesis \eqref{pq} is now replaced by
\bb\label{pq2} \max\{q,s\}<\alpha<p^*.\bbb

We assume that the potential $a\geq 0$ is positive on a subset of $\Omega$ of positive measure and
\bb\label{ahyp} a^{\alpha/(\alpha-s)}\in L^1(\Omega).\bbb

The second main result of this paper establishes an existence and non-existence property if the reaction of problem \eqref{problem1} is dominated by the left-hand side. In this case, solutions exist only for high perturbations of the right-hand side.

\begin{theorem}\label{t2}
Assume that hypotheses \eqref{pq}, \eqref{pq2}, \eqref{ahyp}, and ($\phi_1$)-($\phi_5$) are fulfilled. Then there exists $\Lambda>0$ such that the following properties are true:

(a)  problem \eqref{problem1} does not have any solution for all $0<\lambda<\Lambda$;

(b) problem \eqref{problem1} has a positive solution $U$ for all $\lambda\geq\Lambda$. Moreover,  $U\in C^{1,\mu}(\Omega\cap B_R(0))$ with $\mu =\mu(R)\in(0,1)$.
\end{theorem}

These results remain true if we replace the power-type nonlinearities with general nonlinearities. For instance, the reaction $|u|^{s-2}u$ in the statement of Theorem \ref{t1} corresponding to problem \eqref{problem}, can be replaced by a Carath\'eodory function $f:\RR^N\times\RR\ri\RR$ with the following properties:

($f_1$) $f(x,u)=o(u^{\alpha-1})$ as $u\ri 0^+$, uniformly for a.e. $x\in\RR^N$;

($f_2$) $f(x,u)=O(u^{s-1})$ as $u\ri +\infty$, uniformly for a.e. $x\in\RR^N$.

The above results extend some related properties established by Filippucci, Pucci and R\u adulescu \cite{fpr} in the framework of the $p$-Laplace operator. We also refer to Chorfi and R\u adulescu \cite{chorfi} who studied a related problem driven by the same differential operator and if $\Omega=\RR^N$.

\section{Proof of Theorem \ref{t1}}
We point out that a related property is proved by Azzollini, d'Avenia and Pomponio \cite{azo2} if $\Omega=\RR^N$. However, Theorem~1.3 in \cite{azo2} establishes the existence of a radially symmetric solution and the proof strongly relies on the compact embedding of a Sobolev-type space of functions with radially symmetry into certain Lebesgue spaces. In our case, since $\Omega$ is unbounded but without any symmetry properties, we are not looking for radially symmetric solutions.

The energy functional associated to problem \eqref{problem} is $\calE:\calX\ri\RR$ defined by
$$\calE (u):=\frac 12\int_{\Omega}\phi (|\nabla u|^2)dx+\frac 1\alpha\int_{\Omega}|u|^\alpha dx-\frac{\lambda}{s}\intom |u|^sdx.$$

 By \cite[Proposition 3.1]{azo2}, $\calE$ is well-defined and of class $C^1$. Moreover,  for all $u$, $v\in\calX$ its G\^ateaux directional derivative is given by
$$\langle \calE'(u),v\rangle=\int_{\Omega}\left(\phi'(|\nabla u|^2)\nabla u\nabla v+|u|^{\alpha-2}uv-\lambda |u|^{s-2}uv \right)dx.$$

 We first claim that
\bb\label{claim1}
\mbox{there exists small $r>0$ such that $\inf_{\|u\|=r}\calE (u)>0$.}
\bbb
By ($\phi_1$) we have for all $u\in\calX$
\bb\label{rela1}
\calE (u)\geq\frac{c}{2}\int_{|\nabla u|\leq1}|\nabla u|^qdx+\frac{c}{2}\int_{|\nabla u|>1}|\nabla u|^pdx+\frac 1\alpha\int_{\Omega}|u|^\alpha dx-\frac{\lambda}{s}\intom |u|^sdx.
\bbb
By \cite[Theorem 2.6]{azo2} and our hypothesis $\max\{q,\alpha\}<s<p^*$, it follows that $\calX$ is continuously embedded into $L^s(\Omega)$. So, there exists $c_1>0$ such that
$$\|u\|_s\leq c_1\,\|u\|\quad\mbox{for all}\ u\in\calX.$$
Returning to \eqref{rela1} we obtain for all $u\in\calX$
\bb\label{rela2}
\calE (u)\di\geq \frac{c}{2}\int_{|\nabla u|\leq1}|\nabla u|^qdx+\frac{c}{2}\int_{|\nabla u|>1}|\nabla u|^pdx+\frac 1\alpha\,\|u\|_\alpha^\alpha-\frac{\lambda c_1}{s}\|u\|^s\,.\bbb

Fix $r\in (0,1)$. By \eqref{pq} and \eqref{rela2} we deduce that there are positive constants $c_2,\, c_3$ and $ c_4$ such that for all $u\in\calX$ with $\|u\|=1$
\bb\label{rela3}\begin{array}{ll}
\calE (u)&\di\geq c_2\,(\|\nabla u\|_{p,q}^q+\|u\|_\alpha^\alpha)-c_3\,\|u\|^s\\
&\di\geq c_4\,\left(\|u\|^{\max\{q,\alpha\}}-\|u\|^s \right).\end{array}\bbb
Since $\max\{q,\alpha\}<s<p^*$, relation \eqref{rela3} shows that there exists $c_5>0$ such that
\bb\label{nevo}
\calE (u)\geq c_5\quad\mbox{for all $u\in\calX$ with $\|u\|=r$},\bbb
which proves \eqref{claim1}.

Next, we claim that
\bb\label{claim2}
\liminf_{R\ri +\infty}\calE (u)=-\infty.
\bbb
Indeed, fix $\psi\in\calX\setminus\{0\}$ and $t>0$. Thus, by ($\phi_3$),
$$\begin{array}{ll}
\calE (t\psi)&\di\leq\frac{C}{2}\left(t^q\int_{t|\nabla\psi|\leq 1}|\nabla\psi|^qdx+
t^p\int_{t|\nabla\psi|>1}|\nabla\psi|^pdx\right)+\frac{t^\alpha}{\alpha}\intom |\psi|^\alpha dx-\frac{\lambda t^s}{s}\intom |\psi|^sdx\\
& \ \\
&\di =A_1t^q+A_2t^p+A_3t^\alpha-A_4t^s\ \mbox{(with $A_1, A_2\geq 0$, $A_3,A_4>0$)}\ \ri -\infty\ \mbox{as $t\ri +\infty$},\end{array}$$
by our hypothesis. This proves \eqref{claim2}.

By relations \eqref{claim1}, \eqref{claim2} and using the mountain pass theorem, we find $(u_n)\subset\calX$ such that
\bb\label{unrela}\calE (u_n)\ri c_0\quad\mbox{and}\quad \calE'(u_n)\ri 0\ \mbox{in $\calX^*$}\quad\mbox{as $n\ri\infty$}.\bbb
Here, $c_0:=\inf_{\gamma\in \calC}\max_{t\in [0,1]}\calE (\gamma (t))>0$, where
$$\calC:=\{\gamma:[0,1]\ri\calX;\ \gamma\ \mbox{is continuous},\ \gamma (0)=0,\ \gamma (1)=t_0\psi\},$$
for some fixed $t_0>0$ such that $t_0\|\psi\|>c_5$, where $c_5$ is defined in \eqref{nevo}.

Combining ($\phi_4$) and ($\phi_5$) we deduce that
\bb\label{alphaa}
\phi(t)\leq 2t\phi'(t)\leq s\mu \phi(t)\quad\mbox{for all}\ t>0,
\bbb
hence $s\mu>1$. Thus, by \eqref{alphaa}, $\phi$ is increasing. It follows that
$$\calE(|v|)\leq\calE (v)\quad\mbox{for all}\ v\in\calX.$$
We deduce that we can assume that $u_n\geq 0$ in \eqref{unrela}.

If $(u_n)$ satisfies \eqref{unrela} then
\bb\label{unrela1}\calE(u_n)-\frac{1}{s}\,\langle \calE'(u_n),u_n\rangle =O(1)+o(\|u_n\|)\quad\mbox{as $n\ri\infty$.}\bbb
But, by ($\phi_4$)
\bb\label{unrela2}\begin{array}{ll}
\calE(u_n)-\frac{1}{s}\,\langle \calE'(u_n),u_n\rangle &\di =\intom\left(\frac{1}{2}\, \phi (|\nabla u_n|^2)-\frac{1}{s}\,\phi'(|\nabla u_n|^2)|\nabla u_n|^2 \right)dx\\ &\di+\left(\frac{1}{\alpha}-\frac {1}{s}\right)\intom |u_n|^\alpha dx\\
&\di\geq\frac{1-\mu}{2}\, \intom \phi (|\nabla u_n|^2)dx+\left(\frac{1}{\alpha}-\frac {1}{s}\right)\intom |u_n|^\alpha dx\\
&\di=c_6\intom \phi (|\nabla u_n|^2)dx+c_7\,\|u_n\|^\alpha_\alpha,\end{array}\bbb
where $c_6,\, c_7>0$. 

Next, with an argument similar as the same developed in the first part of this proof, relation \eqref{unrela2} implies that for some $c_8>0$
\bb\label{unrela3}
\calE(u_n)-\frac{1}{s}\,\langle \calE'(u_n),u_n\rangle\geq c_8\,(\|\nabla u_n\|^q_{p,q}+\|u_n\|^\alpha)\quad\mbox{for all}\ n\geq 1.\bbb
Combining \eqref{unrela1} and \eqref{unrela3} we deduce that
$$\|\nabla u_n\|^q_{p,q}+\|u_n\|^\alpha \leq O(1)+o(\|u_n\|)\quad\mbox{as}\ n\ri\infty,$$
which shows that $(u_n)$ is bounded in $\calX$.

Until now we have proved that the Palais-Smale sequence $(u_n)$ of $\calE$ is bounded. Thus, there exists $U\in\calX$ such that, up to a subsequence,
$$u_n\rightharpoonup U\geq 0\quad\mbox{in}\ \calX$$
and
$$u_n\ri U\quad\mbox{in}\ L^s(\Omega)\ \mbox{and}\ L^\alpha(\Omega).$$

We prove in what follows that $U$ is a solution of problem \eqref{problem}. For this purpose we fix $v\in C^\infty_c(\Omega)$ and we set $\omega:={\rm supp}\, v$. Define the functional
$$\calE_0 (u):=\frac 12\int_{\omega}\phi (|\nabla u|^2)dx+\frac 1\alpha\int_{\omega}|u|^\alpha dx.$$
By ($\phi_5$) it follows that $\calE_0$ is convex. Since it is also continuous, it follows that $\calE_0$ is weakly lower semicontinuous. By convexity we have
$$\calE_0(u_n)\leq \calE_0(U)+\langle\calE_0'(u_n),u_n-U\rangle .$$
By \eqref{unrela} we deduce that
$$\limsup_{n\ri\infty}\calE_0(u_n)\leq \calE_0(U).$$
Using now the weakly lower semicontinuity of $\calE_0$ we conclude that
$$\limsup_{n\ri\infty}\calE_0(u_n)= \calE_0(U).$$
Next, with the same arguments as in \cite[p. 210]{azo2}, it follows that
$$\nabla u_n\ri\nabla U\quad\mbox{in}\ L^p(\Omega)+L^q(\Omega).$$
Using \eqref{unrela} and passing to the limit as $n\ri\infty$ we deduce that
$$\int_\omega \phi'(|\nabla U|^2)\nabla U\nabla vdx+\int_\omega |U|^{\alpha-2}Uvdx-\lambda\int_\omega |U|^{s-2}Uvdx=0.$$
By density, this identity is valid for any $v\in\calX$. Thus, $U$ is a solution of problem \eqref{problem}.

We prove in what follows that $U\not=0$. Indeed, if not, it follows that
$$u_n\ri 0\quad\mbox{in}\ L^s(\Omega)\ \mbox{and}\ L^\alpha (\Omega).$$
Thus, by \eqref{unrela}
$$\begin{array}{ll}
\di \frac{c_0}{2}&\di \leq \calE(u_n)-\frac{1}{2}\,\langle \calE'(u_n),u_n\rangle\\
&\di =\frac{1}{2}\intom \left( \phi (|\nabla u_n|^2)-\phi'(|\nabla u_n|^2)|\nabla u_n|^2\right)dx+\left(\frac 1\alpha -\frac 12\right)\intom |u_n|^\alpha dx\\
&\di +\lambda\left(\frac 12-\frac 1s\right)\intom |u_n|^sdx.\end{array}$$
By ($\phi_5$) it follows that $\phi(t^2)-\phi'(t^2)t^2\leq 0$, hence
$$0<\frac{c_0}{2}\leq \left(\frac 1\alpha -\frac 12\right)\intom |u_n|^\alpha dx+\lambda\left(\frac 12-\frac 1s\right)\intom |u_n|^sdx.$$
Passing to the limit as $n\ri\infty$ we get a contradiction. 

We conclude that $U\not=0$ and $U\geq 0$.

\smallskip
(b) By Theorem 1(ii) of Pucci and Servadei \cite{pucserv}, which is based on the Moser iteration, we first deduce that $U\in L^\infty_{loc}(\Omega)$. Next, using the corollary of DiBenedetto \cite[p. 830]{dib}, we conclude that $U\in C^{1,\mu}(\Omega\cap B_R(0))$ with $\mu=\mu(R)\in (0,1)$. A related argument was applied in the proof of Theorem~1 in Yu \cite{yu}.

\smallskip
(c) The function $U\geq 0$ satisfies
$$-\divv [\phi'(|\nabla U|^2)\nabla U]+|U|^{\alpha-2}U\geq 0\ \mbox{in}\ \Omega.$$

We recall that the generalized maximum principle of Pucci and Serrin \cite{pucser04, pucserbook, pucser07} applied to general canonical divergence structure inequalities of the type
$$-{\rm div}\, (A(|\nabla u|)\nabla u)+f(u)\geq 0\quad\mbox{in}\ \Omega,$$
where the function $A=A(t)$ and the nonlinearity $f$ satisfy the following conditions:

\smallskip\noindent
(A1) $A$ is continuous in $\RR^+$;

\smallskip\noindent
(A2) the mapping $t\mapsto tA(t)$ is strictly increasing in $\RR^+$ and $tA(t)\ri 0$ as $t\ri 0^+$;

\smallskip\noindent
(F1) $f\in C(\RR^+_0)$;

\smallskip\noindent
(F2) $f(0)=0$ and $f$ is non-decreasing on some interval $(0,\delta)$, $\delta>0$.

\smallskip
In our case, $f(u)=|u|^{\alpha-2}u$ satisfies (F1) and (F2). We have $A(t)=\phi'(t^2)$ and $t\phi'(t^2)$ is strictly increasing by our hypothesis ($\phi_5$). We also observe that ($\phi_4$) and ($\phi_3$) imply for all $t\in (0,1)$
$$0<t\phi'(t^2)\leq\frac{s\mu}{2}\, \frac{\phi(t^2)}{t}\leq \frac{s\mu}{2}\, t^{q-1}\ri 0\ \mbox{as}\ t\ri 0^+.$$

So, by the Pucci-Serrin maximum principle, we conclude that the non-negative solution $U$ is positive in $\Omega$.
\qed

\section{Proof of Theorem \ref{t2}}
The energy functional associated to problem \eqref{problem1} is $\calF:\calX\ri\RR$ defined by
$$\calF (u):=\frac 12\int_{\Omega}\phi (|\nabla u|^2)dx+\frac 1\alpha\int_{\Omega}|u|^\alpha dx-\frac{\lambda}{s}\intom a(x)\,|u|^sdx.$$

We first establish that $\calF$ is coercive. Indeed, by ($\phi_2$) we have for all $u\in\calX$
\bb\label{calf}\calF (u)\geq \frac{c}{2}\int_{|\nabla u|>1}|\nabla u|^pdx+
\frac{c}{2}\int_{|\nabla u|\leq 1}|\nabla u|^qdx+\frac{1}{\alpha}\,\|u\|^\alpha_\alpha -\frac{\lambda}{s}\intom a(x)|u|^sdx.\bbb

By H\"older's inequality and hypothesis \eqref{ahyp} we obtain
\bb\label{calf1}\intom a(x)|u|^sdx\leq\|a\|_{\alpha/(\alpha-s)}\cdot\left(\intom |u|^\alpha dx \right)^{s/\alpha}=C_1\, \|u\|^s_\alpha,\bbb
where $C_1=C_1(a,\alpha,s,\Omega)$.

Since $\alpha>s$, relations \eqref{calf} and \eqref{calf1} yield
$$\calF (u)\geq \frac{c}{2}\,\|\nabla u\|_{p,q}^p+\frac{1}{\alpha}\,\|u\|^\alpha_\alpha -C_1\, \|u\|^s_\alpha\ri +\infty\quad\mbox{as}\ \|u\|\ri\infty,$$
hence $\calF$ is coercive and bounded from below.

\smallskip
We prove in what follows that problem \eqref{problem1} does not have any solution, provided that $\lambda>0$ is sufficiently small. Indeed, we observe that if $u$ solves \eqref{problem1} then
$$\int_{\Omega}\phi'(|\nabla u|^2)|\nabla u|^2dx+\intom |u|^{\alpha}dx=\lambda\intom a(x)\,|u|^{s}dx.$$

We now estimate the right-hand side of this equality. By hypotheses \eqref{pq2} \eqref{ahyp} we have
$$\begin{array}{ll}
\di \lambda\intom a(x)\,|u|^{s}dx &\di \leq \lambda^{\alpha/(\alpha -s)}\,\frac{\alpha-s}{\alpha}\intom a(x)^{\alpha/(\alpha-s)}dx+\frac{s}{\alpha}\intom |u|^\alpha dx\\
&\di = C(s,a,\alpha)\,\lambda^{\alpha/(\alpha-s)}+\frac{s}{\alpha}\intom |u|^\alpha dx.\end{array}$$

We deduce that if $u$ is a solution of problem \eqref{problem1} then
$$\begin{array}{ll}\di 0&\di \leq \int_{\Omega}\phi'(|\nabla u|^2)|\nabla u|^2dx\leq
 C(s,a,\alpha)\,\lambda^{\alpha/(\alpha-s)}+\left(\frac{s}{\alpha}-1\right)\intom |u|^\alpha dx\\
 &\di <C(s,a,\alpha)\,\lambda^{\alpha/(\alpha-s)},\end{array}$$
 by \eqref{pq2}.

 In conclusion, problem \eqref{problem1} does not have any solution, provided that $\lambda>0$ is small enough. Let
 $$\lambda_*:=\sup\{\lambda>0;\ \mbox{problem \eqref{problem1} does not have a solution}\}>0.$$
 The above arguments show that \eqref{problem1} does not have a solution for all $\lambda<\lambda_*$.
 
In order to obtain sufficient conditions for the existence of solutions, we consider the minimization problem
$$m:=\inf_{u\in\calX}\calF (u)\in\RR.$$
Let $(u_n)\subset\calX$ be a minimizing sequence of $\calF$. Since $\calF(|u_n|)\leq\calF(u_n)$, we can assume that $u_n\geq 0$. Moreover, $(u_n)$ is bounded so, up to a subsequence, we can assume that
$$u_n\rightharpoonup U\geq 0\quad\mbox{in}\ \calX.$$
Hypothesis \eqref{pq2} implies that $\calX$ is compactly embedded into the weighted Lebesgue space $L^s(\Omega;a)$. So, by weak lower semicontinuity and compactness of the embedding, we have
$$\frac 12\int_{\Omega}\phi (|\nabla U|^2)dx+\frac 1\alpha\int_{\Omega}U^\alpha dx\leq\liminf_{n\ri\infty}
\left(\frac 12\int_{\Omega}\phi (|\nabla u_n|^2)dx+\frac 1\alpha\int_{\Omega}u_n^\alpha dx \right)$$
and
$$\intom a(x)\,u_n^sdx\ri \intom a(x)\,U^sdx\quad\mbox{as}\ n\ri\infty.$$

It follows that $U\geq 0$ is a minimizer of $\calF$, that is, $\calF(U)=m$.

We now prove that $U$ is a solution of problem \eqref{problem1}, provided that $\lambda$ is big enough. For this purpose, consider the minimization problem
\bb\label{minpro}
m_0:=\inf_{w\in\calX}\left\{\frac 12\int_{\Omega}\phi (|\nabla w|^2)dx+\frac 1\alpha\int_{\Omega}|w|^\alpha dx;\ 
\frac{1}{s}\intom a(x)|w|^sdx=1 \right\}.\bbb
If $(w_n)\subset\calX$ is a minimizing sequence, then $(w_n)$ is bounded. So, up to a subsequence, we can assume that
$$w_n\rightharpoonup w\quad\mbox{in}\ \calX$$
$$w_n\ri w\quad\mbox{in}\ L^s(\Omega;a).$$
It follows that $w$ is a solution of \eqref{minpro}, hence
$\calF(w)=m_0-\lambda $. We deduce that problem \eqref{problem1} has a solution for all $\lambda>m_0$.

We set
$$\lambda^*:=\inf\{\lambda>0;\ \mbox{problem \eqref{problem1} has a solution}\}.$$
Then $\lambda^*\geq \lambda_*$.

Next, we prove that \eqref{problem1} has a solution for all $\lambda>\lambda^*$. Indeed, if we fix $\lambda>\lambda^*$, then the definition of $\lambda^*$ yields some $\lambda^*<\underline{\lambda}<\lambda$
such that problem \eqref{problem1} has a solution $\underline{U}$ corresponding to $\underline\lambda$. Then $\underline{U}$ is a subsolution of \eqref{problem1}. It remains to prove that problem \eqref{problem1} has a supersolution $\overline{U}$ such that $\overline{U}\geq\underline{U}$. For this purpose we consider the new minimization problem
\bb\label{newpro}
\inf_{v\in\calX}\left\{\frac 12\int_{\Omega}\phi (|\nabla v|^2)dx+\frac 1\alpha\int_{\Omega}|v|^\alpha dx-
\frac{\lambda}{s}\intom a(x)|v|^sdx;\ v\geq\underline{U} \right\}.\bbb

Using the same arguments as above we deduce that the constrained minimization problem \eqref{newpro} has a solution $\overline{U}\geq\underline{U}$.
We conclude that \eqref{problem1} has a solution for all $\lambda>\lambda^*$.

The definition of $\lambda^*$ shows that problem \eqref{problem1} has no solution if $0<\lambda<\lambda^*$. Since $\lambda^*\geq\lambda_*$, we conclude that
$$\lambda^*=\lambda_*=:\Lambda .$$

Until now we know that \eqref{problem1} has no solution if $0<\lambda<\Lambda$ but it has at least one non-negative solution $U$ for all $\lambda\geq\Lambda$. We now prove that problem \eqref{problem1} has a non-negative solution if $\lambda=\Lambda$. Indeed, let $(\lambda_n)$ be a sequence of real numbers such that $\lambda_n\downarrow\Lambda$ as $n\ri\infty$. Let $U_n\geq 0$ be a solution of \eqref{problem1} corresponding to $\lambda_n$. Since $(U_n)\subset\calX$ is bounded, we can assume, passing eventually to a subsequence, that
\bb\label{uu1} U_n\rightharpoonup U_\Lambda\quad\mbox{in}\ \calX\bbb
\bb\label{uu2} U_n\ri U_\Lambda\quad\mbox{in}\ L^s(\Omega;a)\bbb
\bb\label{uu3} U_n\ri U_\Lambda\quad\mbox{a.e.}\ \Omega.\bbb
Since $U_n$ solves \eqref{problem1} for $\lambda=\lambda_n$, it follows that for all $v\in\calX$
\bb\label{uu4}
\int_{\Omega}\phi'(|\nabla U_n|^2)\nabla U_n\nabla vdx+\intom U_n^{\alpha-2}U_nvdx=\lambda_n\intom U_n^{s-2}U_nvdx=0\ \mbox{for all}\ n\geq 1.\bbb
Taking $n\ri\infty$ in \eqref{uu4} and using \eqref{uu1}--\eqref{uu3}, we deduce that $U_\Lambda\geq 0$ is a solution of problem \eqref{problem1} for $\lambda=\Lambda$. We conclude that problem \eqref{problem1} has a solution $U\geq 0$ for every $\lambda\geq\Lambda$. 

Next, as in the proof of Theorem \ref{t1}(b) and using Theorem 1(ii) of Pucci and Servadei \cite{pucserv} in combination with the Moser iteration, we deduce that $U\in L^\infty_{loc}(\Omega)$. This regularity property implies that $U\in C^{1,\mu}(\Omega\cap B_R(0))$, where $\mu =\mu(R)\in(0,1)$
Applying the generalized Pucci-Serrin maximum principle, as in the proof of Theorem \ref{t1}(c), we conclude that $U>0$ in $\Omega$.
\qed

\subsection*{Final comments}  We consider that an interesting research direction with multiple applications concerns the study of nonlinear problems described by the nonlocal term
$$M\left(\int\phi (|\nabla u|^2)|\nabla u|^2 \right),$$
where $\phi$ satisfies hypotheses ($\phi_1$)--($\phi_5$).
Pioneering results have been established by Pucci {\it at all.} \cite{autuori, pucana} in the framework of Kirchhoff problems involving nonlocal operators associated to the standard differential operators.

\medskip
{\bf Acknowledgments.} This research was partially carried out in the ELTE Institutional
Excellence Program (1783-3/2018/FEKUTSRAT) supported by the Hungarian
Ministry of Human Capacities, and it was supported by the Hungarian Scientific Research
Fund OTKA, No. K112157 and SNN125119, and the Slovenian Research Agency Grants
P1-0292, J1-8131, J1-7025, N1-0064, and N1-0083.

\end{document}